\documentclass[12pt]{amsart}
\usepackage{amsmath,amssymb,amsthm,amscd,amsxtra,eucal,latexsym,mathrsfs}

\newtheorem{lem}{Lemma}
\newtheorem{prop}{Proposition}

\newtheorem{thm}{Theorem}

\newtheorem{rem}{Remark}
\theoremstyle{remark}

\DeclareMathOperator{\SL}{SL}

\DeclareMathOperator {\id}{id}

\DeclareMathOperator {\Hom}{Hom}

\renewcommand{\AA}{\mathbb A}

\newcommand {\PP}{\mathbb P}
\newcommand{\QQ}{\mathbb Q}
\newcommand{\BZ}{{\mathbb Z}}
\newcommand{\BN}{{\mathbb N}}

\newcommand{\CE}{{\mathcal E}}

\newcommand{\sk}{{\mathsf k}}

\newcommand{\nc}{\newcommand}
\nc{\on}{\operatorname}
\nc{\ol}{\overline}

\nc{\Fq}{{\mathbb F}_q}
\nc{\Fqb}{\ol{{\mathbb F}_q}}

\nc{\aff}{{\on{aff}}}
\nc{\sph}{{\on{sph}}}
\nc{\Tate}{{\on{Tate}}}
\nc{\Hall}{{\on{Hall}}}
\nc{\bHall}{{\mathbf{Hall}}}
\nc{\CHall}{{{\mathcal H}a\ell\ell}}
\nc{\Mall}{{\on{Mall}}}
\nc{\bMall}{{\mathbf{Mall}}}
\nc{\CMall}{{{\mathcal M}a\ell\ell}}
\nc{\MA}{{\mathcal{MA}}}
\nc{\MC}{{\mathcal{MC}}}
\nc{\Char}{{\mathsf{CH}}}
\nc{\Rep}{{\on{Rep}}}

\nc{\hR}{\hat{\mathcal R}{}^\aff}
\nc{\hRs}{\hat{\mathcal R}{}^\sph}
\nc{\tR}{\tilde{\mathcal R}{}^\sph}
\nc{\tY}{{\tilde{\mathfrak Y}_{n,m}}}
\nc{\fY}{{{\mathfrak Y}_{n,m}}}
\nc{\tX}{{\tilde{\mathfrak X}_{n,m}}}
\nc{\fX}{{{\mathfrak X}_{n,m}}}

\nc{\Gm}{{{\mathbb G}_m}}
\nc{\Pinf}{{{\mathbf P}^{\frac{\infty}{2}}}}
\nc{\Flaff}{{\mathbf{Fl}}}
\nc{\Gr}{{\operatorname{Gr}}}
\nc{\tGr}{\widetilde{\operatorname{Gr}}{}}
\nc{\Pirv}{{{\widetilde{\operatorname{Perv}}}_{L_i,N}}}
\nc{\Parv}{{{\widetilde{\operatorname{Perv}}}_{L_i,\delta_i}}}
\nc{\LB}{{\check{B}}}
\nc{\LG}{{\check{G}}}
\nc{\LL}{{\check{L}}}
\nc{\LT}{{\check{T}}}
\nc{\LR}{{\check{R}}}
\nc{\Lpi}{{\check\Pi}}
\nc{\Lalpha}{{\check\alpha}}
\nc{\Lbeta}{{\check\beta}}
\nc{\Lrho}{{\check\rho}}
\nc{\Gsc}{{G^{\operatorname{sc}}}}
\nc{\Tsc}{{T^{\operatorname{sc}}}}

\nc{\bV}{{\mathbf V}}
\nc{\bVo}{{\overset{\circ}{\mathbf V}}}
\nc{\bP}{{\mathbf P}}
\nc{\bO}{{\mathbf O}}
\nc{\bF}{{\mathbf F}}
\nc{\bI}{{\mathbf I}}
\nc{\bv}{{\mathbf v}}
\nc{\bq}{{\mathbf q}}
\nc{\blambda}{{\boldsymbol{\lambda}}}
\nc{\bmu}{{\boldsymbol{\mu}}}
\nc{\bnu}{{\boldsymbol{\nu}}}
\nc{\bpi}{{\boldsymbol{\pi}}}
\nc{\bfeta}{{\boldsymbol{\eta}}}
\nc{\bzero}{{\boldsymbol{0}}}
\nc{\GO}{{\mathbf{G}_\bO}}
\nc{\GF}{{\mathbf{G}_\bF}}
\nc{\iso}{{\,\widetilde\to\,}}
\nc{\tB}{{\widetilde{\mathcal B}}}

\nc{\BA}{{\mathbb A}}
\nc{\BB}{{\mathbb B}}
\nc{\BC}{{\mathbb C}}
\nc{\BO}{{\mathbb O}}
\nc{\BP}{{\mathbb P}}
\nc{\BQ}{{\mathbb Q}}
\nc{\CA}{{\mathcal A}}
\nc{\CB}{{\mathcal B}}
\nc{\calC}{{\mathcal C}}
\nc{\CY}{{\mathcal Y}}
\nc{\CF}{{\mathcal F}}
\nc{\CG}{{\mathcal G}}
\nc{\CL}{{\mathcal L}}
\nc{\CM}{{\mathcal M}}
\nc{\CN}{{\mathcal N}}
\nc{\CO}{{\mathcal O}}
\nc{\CP}{{\mathcal P}}
\nc{\CS}{{\mathcal S}}
\nc{\CU}{{\mathcal U}}
\nc{\CW}{{\mathcal W}}
\nc{\CX}{{\mathcal X}}

\nc{\fc}{{\mathfrak c}}
\nc{\fg}{{\mathfrak g}}
\nc{\fr}{{\mathfrak r}}
\nc{\fs}{{\mathfrak s}}
\nc{\ft}{{\mathfrak t}}
\nc{\fp}{{\mathfrak p}}
\nc{\fu}{{\mathfrak u}}
\nc{\fv}{{\mathfrak v}}
\nc{\fL}{{\mathfrak L}}

\newcommand{\QED}{$\square$} 
\newenvironment{Prf}{\par\noindent {\it Proof }}{\QED}
\newcommand{\select}[1]{{\it{#1}}}
\newcommand{\RG}{\on{R\Gamma}}

\newcommand{\ZZ}{\BZ}
\newcommand{\Pic}{\on{Pic}}
\nc{\Aut}{\on{Aut}}
\nc{\Spec}{\on{Spec}}
\nc{\wt}{\widetilde}
\nc{\Qlb}{\ol\BQ{}_\ell}
\nc{\getsup}[1]{\stackrel{#1}{\gets}}
\nc{\toup}[1]{\stackrel{#1}{\to}}
\newcommand{\cL}{{\mathcal L}}
\newcommand{\cE}{{\mathcal E}}
\nc{\Perv}{\on{Perv}}
\nc{\PPerv}{\on{{\PP}erv}}
\nc{\Gra}{\on{Gra}}
\nc{\D}{\on{D}}
\renewcommand{\P}{\on{P}}
\nc{\hotimes}{\hat{\otimes}}
\nc{\cF}{\CF}
\nc{\Conv}{\on{Conv}}
\nc{\comp}{\circ}
\newcommand{\hook}[1]{\stackrel{#1}{\hookrightarrow}}
\nc{\ov}[1]{\overline{#1}} 
\nc{\DD}{\mathbb{D}}  %Verdier duality
\nc{\HOM}{{{\mathcal H}om}}
\nc{\R}{\on{R}\!} %highest direct images
\nc{\Vect}{\on{Vect}}
\nc{\GG}{{\mathbb G}}
\nc{\TT}{{\mathbb T}}
\nc{\PSL}{\on{PSL}}
\nc{\Sp}{\on{\mathbb{S}p}}
\nc{\SO}{\on{S\mathbb{O}}}
\nc{\Spin}{\on{\mathbb{S}pin}}

\begin{document}

\title{Twisted geometric Satake equivalence}
\author{Michael Finkelberg and Sergey Lysenko}

\begin{abstract}
Let $\sk$ be an algebraically closed field and 
$\bO=\sk[[t]]\subset \bF=\sk((t))$. For an almost simple algebraic group $G$ 
we classify central extensions $1\to \Gm\to E\to G(\bF)\to 1$; any such 
extension splits canonically over $G(\bO)$. Fix a positive integer $N$ and 
a primitive character $\zeta: \mu_N(\sk)\to\Qlb^*$ (under some assumption on 
the characteristic of $\sk$). Consider the category of $G(\bO)$-biinvariant 
perverse sheaves on $E$ with $\Gm$-monodromy $\zeta$. We show that this is a 
tensor category, which is tensor equivalent to the category of representations 
of a reductive group $\check{G}_{E,N}$. We compute the root datum of 
$\check{G}_{E,N}$.
\end{abstract}
 
\maketitle
\tableofcontents

\section{Introduction}

Let $\sk$ be an algebraically closed field and 
$\bO=\sk[[t]]\subset \bF=\sk((t))$. For an almost simple algebraic group $G$ 
we classify central extensions $1\to \Gm\to E\to G(\bF)\to 1$; any such 
extension splits canonically over $G(\bO)$. Fix a positive integer $N$ and 
a primitive character $\zeta: \mu_N(\sk)\to\Qlb^*$ (under some assumption on 
the characteristic of $\sk$). Consider the category of $G(\bO)$-biinvariant 
perverse sheaves on $E$ with $\Gm$-monodromy $\zeta$. We show that this is a 
tensor category, which is tensor equivalent to the category of representations 
of a reductive group $\check{G}_{E,N}$. We compute the root datum of 
$\check{G}_{E,N}$ in Theorem~\ref{thm_main}. A list of examples is given after Theorem~\ref{thm_main}.

Our result has a natural place in the framework of ``Langlands duality for
quantum groups''~\cite{G}. Namely, if we take $\sk=\BC$, and 
$q=\zeta(\exp(\frac{\pi i}{N}))$ in Conjecture~0.4 of {\em loc. cit.},
then our category of $\zeta$-monodromic perverse sheaves naturally lies inside
the twisted Whittaker sheaves $\on{Whit}^c(\on{Gr}_G)$, 
and corresponds under the equivalence of 
{\em loc. cit.} to the category of representations of the quantum 
Frobenius quotient of $U_q(\check{G})$.
 
From the physical point of view, our result is a 
manifestation of electric-magnetic duality for a {\em rational} parameter
$\Psi$, see~\cite{KW}, Section 11.3. 
Theorem~\ref{thm_main} is a generalization of (\cite{L}, Theorem~3) 
and the classical 
geometric Satake equivalence~\cite{MV}. It was probably known to experts 
for a few years, say it was suggested by an anonymous referee of \cite{L} 
(compare also to \cite{S}). Also, for $G$ simply connected the root data
of ${\check G}_{E,N}$ appeared in~Section~7 of~\cite{Lu}.
Our result should follow essentially by comparing 
Lusztig's results on quantum Frobenius homomorphism on the one hand, and 
Kazhdan-Lusztig-Kashiwara-Tanisaki-Arkhipov-Bezrukavnikov-Ginzburg 
description of representations of quantum groups at roots of unity in 
terms of perverse sheaves on affine Grassmanians, on the other. 
Our goal is 
to provide a short self-contained proof, following the strategy of \cite{MV}.

We are obliged to V.~Drinfeld who explained to us the classification of
central extensions of the loop groups for almost simple groups 
(not necessarily simply connected), see Proposition~\ref{prop_1}.
We are also indebted to R.~Bezrukavnikov and D.~Gaitsgory for useful 
duscussions. M.F. is grateful to University Paris 6 for hospitality and 
support; he was partially
supported by the RFBR grant 09-01-00242 and the Science Foundation of the
SU-HSE awards No.09-08-0008 and 09-09-0009.

\section{Main result}
\label{not}
The general reference for this section is~\cite{BD}~4.5,~5.3 and, for more 
details,~\cite{MV}.

\subsection{Notations} 
\label{sect_notations}
Let $\sk$ be an algebraically closed field of characteristic $p\ge 0$. Let
$G$ be an almost simple algebraic group over $\sk$ with the simply connected cover $\Gsc$.
We fix a Borel subgroup $B\subset G$ and a maximal torus $T\subset B$. We denote the preimage of $T$ in $\Gsc$ by $\Tsc$. The Weyl group of $G,T$ is denoted by $W$. The weight and coweight lattices of $T$ (resp. $\Tsc$) are denoted by $X^*(T)$ and $X_*(T)$ (resp. $X^*(\Tsc)$ and $X_*(\Tsc)$). The root system of $T\subset B\subset G$ is denoted by $R^*\subset
X^*(T)\subset X^*(\Tsc)$; the set of simple roots is
$\Pi^*=\{\Lalpha_1,\ldots,\Lalpha_r\}\subset R^*$. The sum of all positive roots is denoted $2\Lrho$.
The coroot system of $T\subset B\subset G$ is denoted by $R_*\subset X_*(\Tsc)\subset X_*(T)$; the set of simple coroots is
$\Pi_*=\{\alpha_1,\ldots,\alpha_r\}\subset R_*$. Denote by $\fg$ the adjoint representation of $G$.

 There is a unique $W$-invariant bilinear pairing $(,):\
X_*(\Tsc)\times X_*(\Tsc)\to\BZ$ such that $(\alpha,\alpha)=2$ for a
{\em short} coroot $\alpha$. It defines a linear map $\iota:\
X_*(\Tsc)\to X^*(\Tsc)$ such that $(x,y)=\langle x,\iota(y)\rangle$
for any $x,y\in X_*(\Tsc)$. The map $\iota$ extends uniquely by
linearity to the same named map $\iota:\ X_*(T)\to
X^*(T)\otimes_\BZ\BQ$. The bilinear form $(,)$ extends uniquely by linearity to the same named bilinear form $(,):\ X_*(T)\times X_*(T)\to\BQ$. 

\begin{lem}
\label{lem_Coxeter}
 For $\lambda\in X_*(T)$ we have $2\check{h}\iota(\lambda)=\sum_{\check{\alpha}\in R^*}\langle\lambda, \check{\alpha}\rangle\check{\alpha}$, where $\check{h}$ is the dual Coxeter number of $G$. 
\end{lem}
\begin{Prf}
Write  $\Phi_R(.,.)$ for the canonical bilinear $W$-invariant linear form on $X^*(T)$ in the sense of (\cite{Bo}, \textsection 1, section 12). Formula (17) from \select{loc.cit.} reads
$$
2\Phi_R(\check{\beta},\check{\beta})^{-1}=\sum_{\check{\alpha}\in R^+} \langle \check{\alpha}, \beta\rangle^2,
$$
where $\check{\beta}$ is the root corresponding to a short coroot $\beta$. We must check that 
$
\Phi_R(\check{\beta},\check{\beta})=\check{h}^{-1}
$. This is done case by case for all irreducible reduced root systems using the calculation of $\Phi_R$ performed in \textsection 4 of \select{loc.cit.}.  
\end{Prf}

\medskip

 Set $\bO=\sk[[t]],\ \bF=\sk((t))$. The affine Grassmannian $\Gr_G$ is an ind-scheme, the fpqc quotient $G(\bF)/G(\bO)$ (cf. \cite{BD}). Our conventions about $\ZZ/2\ZZ$-gradings are those of \cite{L}. Recall that for free $\bO$-modules of finite type $V_1,V_2$ with an isomorphism $V_1(\bF)\,\iso\, V_2(\bF)$ one has the relative determinant $\det(V_1: V_2)$ (\cite{L}, Section~8.1), this is a $\ZZ/2\ZZ$-graded line given by
$$
\det(V_1: V_2)=\det(V_1/V)\otimes\det(V_2/V)^{-1}
$$
for any $\bO$-lattice $V\subset V_1\cap V_2$. 

 Let $\fL$ be the ($\ZZ/2\ZZ$-graded purely of parity zero) line bundle on $\Gr_G$ whose fibre at $gG(\bO)$ is $\det(\fg(\bO): \fg(\bO)^g)$. Write $\Gra_G$ for the punctured total space (that is, the total space with zero section removed) of the line bundle $\fL$. 
By abuse of notation, the restriction of $\fL$ under the map $G(\bF)\to \Gr_G$, $g\mapsto g G(\bO)$ is again denoted by $\fL$. Write $E^a$ for the punctured total space of the line bundle $\fL$ on $G(\bF)$. Since $\fL$ is naturally a multiplicative $\Gm$-torsor on $G(\bF)$ in the sense of (\cite{DeB}, Section~1.2), $E^a$ gives a central extension 
\begin{equation}
\label{ext_E_a}
1\to\Gm\to E^a\to G(\bF)\to 1,
\end{equation}
it splits canonically over $G(\bO)$. 
   
 Our first result extends the well-known classification of central extensions of $G^{sc}(\bF)$ by $\Gm$ to the almost simple case. % Assume that $p$ does not divide the order of $\pi_1(G)$. This seems not needed!!
 
\begin{prop} 
\label{prop_1}
The isomorphism classes of central extensions of $G(\bF)$ by $\Gm$ are naturally in bijection with integers $m\in\ZZ$ such that $m\iota(X_*(T))\subset X^*(T)$. The corresponding central extension splits canonically over $G(\bO)$, and its group of automorphisms is $\Hom(\pi_0(G(\bF)), \Gm)\,\iso\, \Hom(\pi_1(G), \mu_{\infty}(\sk))$.  
\end{prop}

 Let $d>0$ be the smallest integer such that $d\iota(X_*(T))\subset X^*(T)$, this is a divisor of $2\check{h}$. We pick and denote by $E^c$ the corresponding central extension of $G(\bF)$ by $\Gm$. Any central extension of $G(\bF)$ by $\Gm$ is isomorphic to a multiple of $E^c$. We also pick an isomorphism of central extensions of $G(\bF)$ by $\Gm$ identifying $E^a$ with the $(2\check{h}/d)$-th power of $E^c$. (If $G$ is simply-connected then $d=1$, and the latter isomorphism is uniquely defined). 
 
 Fix a prime $\ell$ different from $p$, write $\P(S)$ (resp., $\D(S)$) for the category of \'etale $\Qlb$-perverse sheaves (resp., the derived category of \'etale $\Qlb$-sheaves) on a $\sk$-scheme (or stack) $S$. Since we are working over an algebraically closed field, we systematically ignore the Tate twists. 

 Fix a positive integer $N$ and a primitive character $\zeta: \mu_N(\sk)\to\Qlb^*$, we assume that $p$ does not divide $2\check{h}N/d$. For the map $s_N: \Gm\to\Gm, x\mapsto x^N$ let $\cL^{\zeta}$ denote the direct summand of $s_{N !}\Qlb$ on which $\mu_N(\sk)$ acts by $\zeta$. If $S$ is a scheme with a $\Gm$-action, 
we say that a perverse sheaf $K$ on $S$ has $\Gm$-monodromy $\zeta$ if it is 
equipped with a $(\Gm, \cL^{\zeta})$-equivariant structure.
  
  Let $\Perv_{G,N}$ denote the category of $G(\bO)$-equivariant $\Qlb$-perverse sheaves on $E^c/G(\bO)$ with $\Gm$-monodromy $\zeta$. 
 
\begin{rem}
\label{Rem_Gm_monodromy}
If $S$ is a scheme and $L$ is a line bundle on $S$, let $\tilde L$ be the total space of the punctured line bundle $L$. Let $f: \tilde L\to \wt{L^m}$ be the map over $S$ sending $l$ to $l^{\otimes m}$. Assume that $p$ does not divide $N$,  
let $\chi: \mu_N(\sk)\to\Qlb^*$ be a character.  Then the functor $K\mapsto f^*K$ is an equivalence between the categories of $\chi$-monodromic perverse sheaves on $\wt{L^m}$ and $\chi^m$-monodromic perverse sheaves on $\tilde L$.
\end{rem}
 
 Pick a primitive character $\zeta_a: \mu_{2\check{h}N/d}(\sk)\to\Qlb^*$ satisfying $\zeta_a^{2\check{h}/d}=\zeta$. % It exists even if $p$ divides $2\check{h}/d$, we do not exlude this case. 
By Remark~\ref{Rem_Gm_monodromy}, $\Perv_{G,N}$ identifies with the category of $G(\bO)$-equivariant perverse sheaves on  $\Gra_G$ with $\Gm$-monodromy $\zeta_a$. 
Set
$$
\PPerv_{G,N}=\Perv_{G,N}[-1]\subset \D(\Gra_G)
$$
For the associativity and commutativity constraints we are going to inroduce later to be natural (and avoid some sign problems), one has to work with $\PPerv_{G,N}$ rather than $\Perv_{G,N}$. Let $\wt\Gr_G$ be the stack quotient of $\Gra_G$ by $\Gm$, where $x\in\Gm$ acts as multiplication by $x^{2\check{h}N/d}$. One may think of $\PPerv_{G,N}$ as the category of certain perverse sheaves on $\wt\Gr_G$. 
 
 We have a natural embedding $X_*(T)\subset\Gr_G$ as the set of
$T$-fixed points. Two coweights $\lambda,\mu\in X_*(T)\subset\Gr_G$ 
lie in the same $G(\bO)$-orbit iff $\lambda\in W(\mu)$. Thus the set of $G(\bO)$-orbits in $\Gr_G$ coincides with the set of Weyl group orbits $X_*(T)/W$, or equivalently, with the set of dominant coweights $X_*^+(T)\subset X_*(T)$. The orbit corresponding to $\lambda\in X_*^+(T)$ will be denoted by $\Gr^\lambda$. The $G$-orbit of $\lambda$ is isomorphic to a partial flag variety $\CB^\lambda=G/P^\lambda$
where $P^\lambda$ is a parabolic subgroup whose Levi has the Weyl group $W^\lambda\subset W$ coinciding with the stabilizer of $\lambda$
in $W$. Write $\Gra_G^{\lambda}$ for the preimage of $\Gr_G^{\lambda}$ under $\Gra_G\to \Gr_G$.

  Let $\Aut(\bO)$ denote the group ind-scheme over $\sk$ such that, for any $\sk$-algebra $R$, $\Aut(\bO)(R)$ is the automorphism group of the topological $R$-algebra $R\hotimes\bO$ (\cite{BD}, 2.6.5). Write $\Aut^0(\bO)$ for the reduced part of $\Aut(\bO)$. 
The group scheme $\Aut^0(\bO)$ acts naturally on the exact sequence (\ref{ext_E_a}) acting trivially on $\Gm$ and preserving $G(\bO)$.  

The action of the loop rotation group
$\Gm\subset\Aut^0(\bO)$ contracts $\Gr^\lambda$ to
$\CB^\lambda\subset\Gr_\lambda$, and realizes $\Gr^\lambda$ as a composition of affine fibrations over $\CB^\lambda$. We denote the projection $\Gr^\lambda\to\CB^\lambda$ by $\varpi_\lambda$. 

 If $\check{\nu}\in X^*(\Tsc)$ is orthogonal to all coroots $\alpha$ satisfying $\langle\lambda, \check{\alpha}\rangle=0$ then write $\CO(\check{\nu})$ for the corresponding $\Gsc$-equivariant line bundle on $\CB^{\lambda}$. It is canonically trivialized at $1\in\CB^{\lambda}$. If $\lambda\in X_*(T)$ then $\iota(\lambda)\in X^*(\Tsc)$ gives
rise to the line bundle $\CO(\iota(\lambda))$ on $\CB^\lambda$.  

 For a free $\bO$-module $\cE$ write $\cE_{\bar c}$ for its geometric fibre. Write $\Omega$ for the completed module of relative differentials of $\bO$ over $\sk$.
 
  For a root $\check{\alpha}$ write $\fg^{\check{\alpha}}\subset \fg$ for the corresponding root subspace. If $\check{\alpha}=\sum_{i=1}^r a_i\check{\alpha}_i$ we have $\fg^{\check{\alpha}}\iso
\otimes_{i=1}^r (\fg^{\check{\alpha}_i})^{a_i}$ canonically. In particular, $\fg^{-\check{\alpha}}$ is identified with the dual of $\fg^{\check{\alpha}}$ via the Killing form. 

 Pick a trivialization $\phi_i: \fg^{\check{\alpha}_i}\,\iso\, \sk$ for each simple root $\check{\alpha}_i\in\Pi^*$. Set $\Phi=\{\phi_i\}_{i=1}^r$. 
 
\begin{lem} 
\label{lem_fibres_fL}
Let $\lambda\in X_*(T)$. The family $\Phi$ yields a uniquely defined $\ZZ/2\ZZ$-graded 
isomorphism 
$$
\fL\mid_{\Gr^{\lambda}_G}\simeq \Omega_{\bar c}^{\check{h}\langle \lambda, \iota(\lambda)\rangle}\otimes
\varpi_\lambda^* \CO(2\check{h} \iota(\lambda))
$$
\end{lem}

 Set $X^{*+}(\check{T}_N)=\{\lambda\in X_*^+(T)\mid d\iota(\lambda)\in N X^*(T)\}$. By Lemma~\ref{lem_fibres_fL}, for $\lambda\in X_*^+(T)$ the scheme $\Gra^{\lambda}_G$ admits a $G(\bO)$-equivariant local system with $\Gm$-monodromy $\zeta_a$ iff $\lambda\in X^{*+}(\check{T}_N)$. 
 
\begin{rem}  In (\cite{BD}, Section~4.4.9, p. 166, formula (217)) an extension of $G(\bF)$ by $\Gm$ has been constructed whose square identifies with (\ref{ext_E_a}). So, $d$ is a divisor of $\check{h}$. Another way to see this is to note that, by Lemma~\ref{lem_Coxeter}, for $\lambda\in X_*(T)$ we have $\check{h}\iota(\lambda)=\sum_{\check{\alpha}\in R^{*+}} \langle \lambda, \check{\alpha}\rangle\check{\alpha}
\in X^*(T)$. Here $R^{*+}$ denotes the set of positive roots.
\end{rem} 

 Write $\Omega^{\frac{1}{2}}(\bO)$ for the groupoid of square roots of $\Omega$. For $\CE\in \Omega^{\frac{1}{2}}(\bO)$ and $\lambda\in X^{*+}(\check{T}_N)$ 
define the line bundle $\cL_{\lambda, \CE}$ on $\Gr^{\lambda}_G$ as 
$$
\cL_{\lambda,\CE}=\CE_{\bar c}^{\frac{d}{N}(\lambda, \lambda)}\otimes \varpi_\lambda^* \CO(\frac{d}{N}\iota(\lambda))
$$
It is equipped with an isomorphism $\cL_{\lambda,\CE}^{2\check{h}N/d}\,\iso\, \fL\mid_{\Gr^{\lambda}_G}$. Write ${\overset{\circ}{\cL}}_{\lambda,\CE}$ for the punctured total space of $\cL_{\lambda,\CE}$. Denote by 
\begin{equation}
\label{map_p_lambda}
p_{\lambda}: {\overset{\circ}{\cL}}_{\lambda,\CE}\to \Gra_G^{\lambda}
\end{equation}
the map over $\Gr_G^{\lambda}$ sending $x$ to $x^{2\check{h}N/d}$. 

 For $\lambda\in X^{*+}(\check{T}_N)$ we define $\CA^{\lambda}_{\CE}\in\PPerv_{G,N}$ as the intermediate extension of $E^{\lambda}_{\CE}[-1+\dim\Gra^{\lambda}_G]$. 
Here $E^{\lambda}_{\CE}$ is the local system on $\Gra_G^{\lambda}$ with $\Gm$-monodromy $\zeta_a$ equipped with an isomorphism
$p_{\lambda}^* E^{\lambda}_{\CE}\,\iso\, \Qlb$. Both $E^{\lambda}_{\CE}$ and $\CA^{\lambda}_{\CE}$ are defined up to a unique isomorphism. The irreducible objects of $\PPerv_{G,N}$ are exactly $\CA^{\lambda}_{\CE}$, $\lambda\in X^{*+}(\check{T}_N)$.
  
  As in (\cite{L}, Proposition~11) one shows that each $\CA^{\lambda}_{\CE}$ has nontrivial usual cohomology sheaves only in degrees of the same parity, and derives from this that $\PPerv_{G,N}$ is semi-simple.
   
\subsection{Convolution} 
\label{sect_convolution}
Consider the automorphism $\tau$ of $E^a\times E^a$ sending $(g,h)$ to $(g, gh)$. Let $G(\bO)\times G(\bO)\times\Gm$ act on $E^a\times E^a$ in such a way that $(\alpha,\beta, b)\in G(\bO)\times G(\bO)\times\Gm$ send $(g,h)$ to $(g\beta^{-1}b^{-1}, \beta bh\alpha)$. Write $E^a\times_{G(\bO)\times\Gm} \Gra_G$ for the quotient of $E^a\times E^a$ under this free action. Then $\tau$ induces an isomorphism
$$
\bar\tau: E^a\times_{G(\bO)\times\Gm} \Gra_G\,\iso\, \Gr_G\times\Gra_G
$$
sending $(g, hG(\bO))$ to $(\bar gG(\bO), ghG(\bO))$, where $\bar g$ is the image of $g\in E^a$ in $G(F)$. 
Set $m$ be the composition of $\bar\tau$ with the projection to $\Gra_G$. Let $p_G: E^a\to \Gra_G$ be the map $h\mapsto hG(\bO)$. Similarly to (\cite{MV}, \cite{L}), we get a diagram
$$
\Gra_G\times\Gra_G\getsup{p_G\times\id} E^a\times \Gra_G \toup{q_G} E^a\times_{G(\bO)\times\Gm} \Gra_G \toup{m} \Gra_G,
$$
where $q_G$ is the quotient map under the action of $G(\bO)\times\Gm$.  
 
 For $K_1,K_2\in \PPerv_{G,N}$ define the convolution product $K_1\ast K_2\in \D(\Gra_G)$ by $K_1\ast K_2=m_! K$, where $K[1]$ is a perverse sheaf on $E^a\times_{G(\bO)\times\Gm} \Gra_G$ equipped with an isomorphism $q_G^*K\,\iso\, p_G^*K_1\boxtimes K_2$. Since $q_G$ is a $G(\bO)\times\Gm$-torsor and $p_G^*K_1\boxtimes K_2$ is naturally equivariant under $G(\bO)\times\Gm$, $K$ is defined up to a unique isomorphism. 
 
\begin{lem} If $K_1,K_2\in \PPerv_{G,N}$ then $K_1\ast K_2\in \PPerv_{G,N}$.
\end{lem}
\begin{Prf} Following \cite{MV}, stratify $E^a\times_{G(\bO)\times\Gm}\Gra_G$ by locally closed subschemes $p_G^{-1}(\Gra_G^{\lambda})\times_{G(\bO)\times\Gm} \Gra_G^{\mu}$ for $\lambda,\mu\in X_*^+(T)$. Stratify $\Gra_G$ by $\Gra_G^{\lambda}$, $\lambda\in X_*^+(T)$. By (\cite{MV}, Lemma~4.4), $m$ is a stratified semi-small map, our assertion follows.
\end{Prf} 
 
\medskip

 In a similar way one defines a convolution product $K_1\ast K_2\ast K_3$ of $K_i\in \PPerv_{G,N}$. Moreover, $(K_1\ast K_2)\ast K_3\,\iso\, K_1\ast K_2\ast K_3\,\iso\, K_1\ast (K_2\ast K_3)$ canonically, and   $\CA^0_{\CE}$ is a unit object in $\PPerv_{G,N}$.  
 
\subsection{Fusion} 
\label{sect_Fusion} As in \cite{MV}, we will show that the convolution product on $\PPerv_{G,N}$ can be interpreted as a fusion product, thus leading to a commutativity constraint on $\PPerv_{G,N}$.

 Fix $\CE\in\Omega^{\frac{1}{2}}(\bO)$ and consider the group scheme $\Aut_2(\bO):=\Aut(\bO, \CE)$ as in (\cite{BD}, 3.5.2). It fits into an exact sequence $1\to \mu_2\to \Aut_2(\bO)\to \Aut(\bO)\to 1$, and $\Aut_2(\bO)$ is connected. Write $\Aut^0_2(\bO)$ for the preimage of $\Aut^0(\bO)$ in $\Aut_2(\bO)$. 

 The map (\ref{map_p_lambda}) is $\Aut_2^0(\bO)$-equivariant, so the action of $\Aut^0(\bO)$ on $\Gra_G$ lifts to a $\Aut_2^0(\bO)$-equivariant structure on each $\CA^{\lambda}_{\CE}\in\PPerv_{G,N}$. The corresponding  $\Aut_2^0(\bO)$-equivariant structure on each $\CA^{\lambda}_{\CE}$ is unique, as the action of $\Aut_2^0(\bO)$ on $\ov{\Gra}^{\lambda}_G$ factors through a smooth connected quotient group of finite type. Here $\ov{\Gra}^{\lambda}_G$ is the preimage of $\ov{\Gr}^{\lambda}_G$ under the projection $\Gra_G\to \Gr_G$. 
 
 Let $X$ be a smooth connected projective curve over $k$. For a closed $x\in X$ let $\bO_x$ be the completed local ring of $X$ at $x$, and $\bF_x$ its fraction field. Write $\cF_G^0$ for the trivial $G$-torsor on a scheme (or stack). Write $\Gr_{G,x}=G(\bF_x)/G(\bO_x)$ for the corresponding affine grassmanian. Then $\Gr_{G,x}$ identifies canonically with the ind-scheme classifying a $G$-torsor $\cF_G$ on $X$ together with a trivialization $\nu: \cF_G\,\iso\, \cF^0_G\mid_{X-x}$. 
 
 For $m\ge 1$ write $\Gr_{G, X^m}$ for the ind-scheme classifying $(x_1,\ldots, x_m)\in X^m$, a $G$-torsor $\cF_G$ on $X$, and a trivialization $\cF_G\,\iso\, \cF_G^0\mid_{X-\cup x_i}$. 
 
 Let $G_{X^m}$ be the group scheme over $X^m$ classifying $\{(x_1,\ldots, x_m)\in X^m, \mu\}$, where $\mu$ is an automorphism of $\cF_G^0$ restricted to the formal neighborhood of $D=x_1\cup\ldots\cup x_m$ in $X$. The fibre of $G_{X^m}$ over $(x_1,\ldots, x_m)\in X^m$ is $\prod_i G(\bO_{y_i})$ with $\{y_1,\ldots, y_s\}=\{x_1,\ldots, x_m\}$ and $y_i$ pairwise distinct.
 
 Let $\fL_{X^m}$ be the ($\ZZ/2\ZZ$-graded purely of parity zero) line bundle on $\Gr_{G, X^m}$ whose fibre is $\det\RG(X, \fg_{\cF^0_G})\otimes\det\RG(X, \fg_{\cF_G})^{-1}$. Here for a $G$-module $V$ and a $G$-torsor $\cF_G$ on a base $S$ we write $V_{\cF_G}$ for the induced vector bundle on $S$. 
 
\begin{lem} For a $\sk$-point $(x_1,\ldots, x_m,\cF_G)$ of $\Gr_{G, X^m}$ let $\{y_1,\ldots, y_s\}=\{x_1,\ldots, x_m\}$ with $y_i$ pairwise distinct. The fibre of $\fL_{X^m}$ at this $\sk$-point is canonically isomorphic as $\ZZ/2\ZZ$-graded to
$$
\otimes_{i=1}^s \det(\fg(\bO_{y_i}): \fg_{\cF_G}(\bO_{y_i}))
$$
\end{lem}  
 
 Write $\Gra_{G, X^m}$ for the punctured total space of $\fL_{X^m}$. The group scheme $G_{X^m}$ acts naturally on $\Gra_{G, X^m}$ and $\Gr_{G, X^m}$, and the projection $\Gra_{G, X^m}\to \Gr_{G, X^m}$ is $G_{X^m}$-equivariant. Let $\Perv_{G, N, X^m}$ be the category of $G_{X^m}$-equivariant perverse sheaves on $\Gra_{G, X^m}$ with $\Gm$-monodromy $\zeta_a$. Set
$$
\PPerv_{G,N,X^m}=\Perv_{G,N,X^m}[-m-1]\subset \D(\Gra_{G, X^m})
$$ 
 
 For $x\in X$ write $D_x=\Spec\bO_x$, $D_x^*=\Spec\bF_x$. Consider the diagram, where the left and right square is cartesian
$$
\begin{array}{ccccccc}
\Gra_{G,X}\times\Gra_{G,X} & \getsup{\tilde p_{G,X}} &\tilde C_{G,X} &  \toup{\tilde q_{G,X}} & \wt\Conv_{G,X} & \toup{\tilde m_X} & \Gra_{G, X^2}\\
\downarrow && \downarrow && \downarrow && \downarrow\\
\Gr_{G,X}\times\Gr_{G,X} & \getsup{p_{G,X}} & C_{G,X} &  \toup{q_{G,X}} & \Conv_{G,X} & \toup{m_X} & \Gr_{G, X^2}
\end{array}
$$
Here the low row is the convolution diagram from \cite{MV}. Namely, $C_{G,X}$ is the ind-scheme classifying 
collections:
\begin{equation}
\label{point_of_C_GX}
\left\{ 
\begin{array}{l}
x_1, x_2\in X, \; G-\mbox{torsors} \; \cF^1_G, \cF^2_G \;\mbox{on} \; X \; \mbox{with} \; \nu_i: \cF^i_G\,\iso\, \cF^0_G\mid_{X-x_i},\\
\mu_1: \cF^1_G\,\iso\, \cF^0_G\mid_{D_{x_2}}
\end{array}
\right.
\end{equation} 
The map $p_{G,X}$ forgets $\mu_1$. The ind-scheme $\Conv_{G,X}$ classifies collections:
\begin{equation}
\label{point_of_Conv}  
\left\{ 
\begin{array}{l}
x_1, x_2\in X, \; G-\mbox{torsors} \; \cF^1_G, \cF_G \;\mbox{on} \; X,  \\
\mbox{isomorphisms}\; \nu_1: \cF^1_G\,\iso\, \cF^0_G\mid_{X-x_1}\; \mbox{and}\; \eta: \cF^1_G\,\iso\, \cF_G\mid_{X-x_2}
\end{array}
\right.
\end{equation}
The map $m_X$ sends this collection to $(x_1, x_2, \cF_G)$ together with the trivialization $\eta\comp\nu_1^{-1}: \cF^0_G\,\iso\, \cF_G\mid_{X-x_1-x_2}$. 

 The map $q_{G,X}$ sends (\ref{point_of_C_GX}) to (\ref{point_of_Conv}), where $\cF_G$ is obtained by gluing $\cF^1_G$ on $X-x_2$ and $\cF_G^2$ on $D_{x_2}$ using their identification over $D_{x_2}^*$ via $\nu_2^{-1}\comp\mu_1$.   
 
  The canonical $\ZZ/2\ZZ$-graded isomorphism 
$$
q_{G,X}^*m_X^*\fL_{X^2}\,\iso\, p_{G,X}^*(\fL_X\boxtimes\fL_X)
$$
allows to define $\tilde q_{G,X}$, it sends (\ref{point_of_C_GX}) together with $v_i\in\det(\fg(\bO_{x_i}): \fg_{\cF_G^i}(\bO_{x_i}))$ for $i=1,2$ to the image of (\ref{point_of_C_GX}) under $q_{G,X}$ together with $v_1\otimes v_2$. Here we used the isomorphism
\begin{multline}
\label{fibre_fL_for_fusion}
\det(\fg(\bO_{x_1}): \fg_{\cF_G^1}(\bO_{x_1}))\otimes \det(\fg(\bO_{x_2}): \fg_{\cF_G^2}(\bO_{x_2}))\,\iso\\ 
\det(\fg(\bO_{x_1}): \fg_{\cF_G^1}(\bO_{x_1})\otimes
\det(\fg_{\cF_G^1}(\bO_{x_2}): \fg_{\cF_G}(\bO_{x_2}))
\end{multline}
given by $\mu_1$ and $\fg_{\cF_G}(\bO_{x_2})\,\iso\, \fg_{\cF^2_G}(\bO_{x_2})$, so (\ref{fibre_fL_for_fusion}) is the fibre of $\fL_{X^2}$ over $\cF_G$. 
 
 For $K_1,K_2\in \PPerv_{G,N,X}$ there is a (defined up to a unique isomorphism) perverse sheaf $K_{12}[3]$ on $\wt\Conv_{G,X}$ with $\tilde q_{G,X}^*K_{12}\,\iso\, \tilde p_{G,X}^*(K_1\boxtimes K_2)$. Moreover, $K_{12}$ has $\Gm$-monodromy $\zeta_a$. We then let
$$
K_1\ast_X K_2= \tilde m_{X !} K_{12}
$$
Let $U\subset X^2$ be the complement to the diagonal. Let $j: \Gra_{G, X^2}(U)\hook{} \Gra_{G,X^2}$ be the preimage of $U$. Recall that $m_X$ is stratified small, an isomorphism over the preimage of $U$ (\cite{MV}), so the same holds for $\tilde m_X$. Thus, $(K_1\ast_X K_2)[3]$ is a perverse sheaf, the intermediate extension from $\Gra_{G, X^2}(U)$. Clearly, $K_1\ast_X K_2$ is $G_{X^2}$-equivariant over $\Gra_{G, X^2}(U)$, and this property is preserved under the intermediate extension. So, $K_1\ast_X K_2\in \Perv_{G,N, X^2}$. 

  Let $\Omega_X$ be the canonical line bundle on $X$. Write $\Omega_X^{\frac{1}{2}}(X)$ for the groupoid of square roots of $\Omega_X$. For $\CE_X\in \Omega_X^{\frac{1}{2}}(X)$ let $\hat X_2\to X$ be the $\Aut_2^0(\bO)$-torsor whose fibre over $x$ is the scheme of isomorphisms between $(\CE_x, \bO_x)$ and $(\CE, \bO)$. Then $\Gr_{G,X}\,\iso\, \hat X_2\times_{\Aut_2^0(\bO)} \Gr_G$ (cf. \cite{BD}, 5.3.11), and similarly $\Gra_{G,X}\,\iso\, \hat X_2\times_{\Aut_2^0(\bO)} \Gra_G$. Since any $K\in\Perv_{G,N}$ is $\Aut^0_2(\bO)$-equivariant (in a unique way), we get a fully faithful functor 
\begin{equation}
\label{functor_tau_0}
\tau^0: \PPerv_{G,N}\to \PPerv_{G,N,X}
\end{equation}
sending $K$ to the descent of $\Qlb\boxtimes K$ under $\hat X_2\times \Gra_G\to \Gra_{G,X}$.

   Let $i: \Gra_{G,X}\to \Gra_{G,X^2}$ be the preimage of the diagonal in $X^2$. For $F_i\in \PPerv_{G,N}$ letting $K_i=\tau^0 F_i$ define 
$$
K_{12}\mid_U:=K_{12}\mid_{\Gra_{G,X^2}(U)}
$$ 
as above (it is placed in perverse degree 3). We get $K_1\ast_X K_2\,\iso\, j_{!*}(K_{12}\mid_U)$ and 
$
\tau^0(F_1\ast F_2)\,\iso\, i^*(K_1\ast_X K_2)
$. So, the involution $\sigma$ of $\Gra_{G,X^2}$ interchanging $x_i$ yields   
$$
\tau^0(F_1\ast F_2)\,\iso\,i^*j_{!*}(K_{12}\mid_U)\,\iso\, i^*j_{!*}(K_{21}\mid_U)\,\iso\, \tau^0(F_2\ast F_1),
$$
because $\sigma^*(K_{12}\mid_U)\,\iso\, K_{21}\mid_U$ canonically.  As in (\cite{BD}, 5.3.13-5.3.17) one shows that the associativity and commutativity constraints are compatible. Thus, $\PPerv_{G,N}$ is a symmetric monoidal category.   

 The idea to use $\tau^0$ instead of $\tau^0[1]$ in the above definition of the commutativity constraint goes back to (\cite{BD}, 5.3.17), this is a way to avoid sign problems.
   
\begin{rem} 
\label{rem_rigidity}
Write $\P_{G(\bO)}(\Gra_G)$ for the category of $G(\bO)$-equivariant perverse sheaves on $\Gra_G$. Let $\star$ be the covariant self-functor on $\P_{G(\bO)}(\Gra_G)$ induced by the map $E^a\to E^a$, $e\mapsto e^{-1}$. Then $K\mapsto K^{\vee}:=\DD(\star K)[-2]$ is a contravariant functor $\PPerv_{G,N}\to\PPerv_{G,N}$. As in (\cite{L}, Remark~6), one checks that for $K_i\in\PPerv_{G,N}$ we have canonically $\R\Hom(K_1\ast K_2, K_3)\,\iso\, \R\Hom(K_1, K_3\ast K_2^{\vee})$. So, $K_3\ast K_2^{\vee}$ represents the internal $\HOM(K_2, K_3)$ in the sense of the tensor structure on $\PPerv_{G,N}$. Besides, $\star(K_1\ast K_2)\,\iso\, (\star K_2)\ast (\star K_1)$ canonically.
\end{rem} 

\subsection{Main result}
\label{sect_main_result}
 In Section~\ref{subsection_fibre_functor} below we introduce a tensor category $\PPerv^{\natural}_{G,N}$ obtained from $\PPerv_{G,N}$ by some modification of the commutativity constraint. Set 
$$
X^*(\check{T}_N)=\{\nu\in X_*(T)\mid d\iota(\nu)\in N X^*(T)\}$$ 
Let $\check{T}_N=\Spec k[X^*(\check{T}_N)]$ be the torus whose weight lattice is $X^*(\check{T}_N)$. The natural inclusion $X^*(T)\subset X_*(\check{T}_N)$ allows to see each root $\check{\alpha}\in R^*$ as a coweight of $\check{T}_N$. For $a\in\QQ, a>0$ written as $a=a_1/a_2$ with $a_i\in\BN$ prime to each other say that $a_2$ is \select{the denominator} of $a$. Recall that $p$ does not divide $2\check{h}N/d$.

\begin{thm} 
\label{thm_main}
There is a connected semi-simple group $\check{G}_N$ and a canonical equivalence of tensor categories
$$
\PPerv^{\natural}_{G,N}\,\iso\, \Rep(\check{G}_N)
$$
There is a canonical inclusion $\check{T}_N\subset\check{G}_N$ whose image is a maximal torus in $\check{G}_N$. The Weyl groups of $G$ and $\check{G}_N$ 
viewed as subgroups of $\Aut(X^*(\check{T}_N))$ are the same. Our choice of a Borel subgroup $T\subset B\subset G$ yields a Borel subgroup $\check{T}_N\subset \check{B}_N\subset\check{G}_N$. The corresponding simple roots (resp., coroots) of $(\check{G}_N,\check{T}_N)$ are $\delta_i\alpha_i$ (resp., $\frac{\check{\alpha}_i}{\delta_i}$), where $\delta_i$ is the denominator of $\frac{d(\alpha_i,\alpha_i)}{2N}$.
\end{thm}

\noindent
{\bf Examples.} (If $G$ is simply-connected then $d=1$.)
\begin{itemize}
\item $G=\SL_2$ then $\check{G}_N\,\iso\, \SL_2$ for $N$ even, 
and $\check{G}_N\,\iso\, \PSL_2$ for $N$ odd. 
\item $G=\PSL_2$ then $d=2$, and 
$\check{G}_N\,\iso\, \SL_2$ for $N$ odd, $\check{G}_N\,\iso\, \PSL_2$ 
for $N$ even. 
\item $G=\Sp_{2n}$ then $\check{G}_N\,\iso\, \SO_{2n+1}$ for $N$ odd, and  $\check{G}_N\,\iso\,\Sp_{2n}$ for $N$ even. For $N=2$ this has been also proved in \cite{L}.
\item $G=\Spin_{2n+1}$ with $n\ge 2$ then 
$$
\check{G}_N\,\iso\,\left\{
\begin{array}{lc}
\Sp_{2n}/\{\pm 1\}, & N \;\mbox{odd}\\
\Spin_{2n+1}, & N\; \mbox{even and}\; nN/2\; \mbox{even}\\
\SO_{2n+1}, & N\; \mbox{even and}\; nN/2\; \mbox{odd}
\end{array}
\right.
$$
\item $G=G_2$ has trivial center, and $\check{G}_N\,\iso\, G_2$ for any $N$
\item $G=F_4$ has trivial center, and $\check{G}_N\,\iso\, F_4$ for any $N$
\item $G=E_8$ has trivial center, and $\check{G}_N\,\iso\, E_8$ for any $N$
\item $G$ simply-connected of type $E_6$, its center identifies with $\ZZ/3\ZZ$ and 
$$
\check{G}_N\,\iso\,\left\{
\begin{array}{lc}
\mbox{adjoint of type}\; E_6, & 3\nmid N\\
\mbox{simply-connected of type}\; E_6, & 3\mid N
\end{array}
\right.
$$
\item $G$ is simply-connected of type $E_7$, its center identifies with $\ZZ/2\ZZ$ and 
$$
\check{G}_N\,\iso\,\left\{
\begin{array}{lc}
\mbox{simply-connected of type}\; E_7, & N \;\mbox{even}\\
\mbox{adjoint of type}\; E_7, & N\; \mbox{odd}
\end{array}
\right.
$$
\end{itemize} 
 
\begin{rem} The case when $p$ divides $2\check{h}/d$, but $p$ does not divide $N$ can also be treated. In this case one can pick a character $\zeta_a: \mu_{2\check{h}N/d}(\sk)\to \Qlb^*$ satisfying $\zeta_a^{2\check{h}/d}=\zeta$ and define $\PPerv_{G,N}^{\natural}$ in the same way. 
But the category  $\PPerv_{T,G,N}$ will have more objects, we excluded this case to simplify the proof.
\end{rem}

\medskip
 
\section{Classification of central extensions}
\label{section_extensions}

\subsection{Simply-connected case} In this subsection we remind the classification of central extensions of $G^{sc}(\bF)$ by $\Gm$ in relation with \cite{DeB}. 

 By \cite{DeB}, the central extensions of $G^{sc}$ by (the sheaf version of) $K_2$ are classified by integer-valued $W$-invariant quadratic forms on $X_*(T^{sc})$ (and have no automorphisms). Let $Q$ be the unique $\ZZ$-valued quadratic form on $X_*(T^{sc})$ satisfying $Q(\alpha)=1$ for a short coroot $\alpha$. So, $(\lambda_1,\lambda_2)=Q(\lambda_1+\lambda_2)-Q(\lambda_1)-Q(\lambda_2)$ for $\lambda_i\in X_*(T^{sc})$. Let 
\begin{equation}
\label{ext_can_K2}
1\to K_2\to E_Q\to G^{sc}\to 1
\end{equation}
denote the central extension corresponding to $Q$. 

 Write $v(f)$ for the valuation of $f\in\bF^*$. Write $(.,.)_{st}$ for the tame symbol given by
$$
(f,g)_{st}=(-1)^{v(f)v(g)}(g^{v(f)}f^{-v(g)})(0)
$$ 
for $f,g\in F^*$. We may view it as a map $K_2(\bF)\to \sk^*$. Taking the $\bF$-valued points of (\ref{ext_can_K2}) and further the push-forward by the tame symbol, one gets a central extension
\begin{equation}
\label{ext_sc_case_k_points}
1\to \sk^*\to \bar G\to G^{sc}(\bF)\to 1
\end{equation}
 
 For $\theta\in\pi_1(G)$ write $\Gr_G^{\theta}$ for the connected component of $\Gr_G$ that contains $t^{\lambda}G(\bO)$ for $\lambda\in X_*(T)$ whose image in $\pi_1(G)$ equals $\theta$. The natural map $\Gr_{G^{sc}}\to\Gr_G^0$ is an isomorphism. 
 
 From \cite{F} one knows that there is a line bundle $\cL$ on $\Gr_G$ generating the Picard group $\Pic(\Gr^{\theta}_G)\,\iso\, \ZZ$ of each connected component $\Gr_G^{\theta}$ of $\Gr_G$, and an isomorphism $\cL^{2\check{h}}\,\iso\, \fL$. 
 Write $\bar G_Q$ for the punctured total space of the inverse image of $\cL$ under $G^{sc}(F)\to\Gr_G$, $x\mapsto xG(\bO)$. It can be seen as the Mumford extension 
\begin{equation} 
\label{ext_Mumford_sc_case}
1\to \Gm\to \bar G_Q\to G^{sc}(\bF)\to 1,
\end{equation}
that is, the ind-scheme classifying pairs $(g\in G^{sc}(\bF), \sigma)$, where $\sigma: g^*\cL\,\iso\, \cL$ is an isomorphism over $\Gr_G^0$. The central extension (\ref{ext_Mumford_sc_case}) splits canonically over $G^{sc}(\bO)$. Any central extension of $G^{sc}(\bF)$ by $\Gm$ is a multiple of (\ref{ext_Mumford_sc_case}) and has no automorphisms. In the rest of this section we prove the following.

\begin{lem} Passing to $\sk$-points in (\ref{ext_Mumford_sc_case}) one gets a central extension isomorphic to (\ref{ext_sc_case_k_points}).
\end{lem}

 For a central extension $1\to A\to E\to H\to 1$ we write $(.,.)_c: H\times H\to A$ for the corresponding commutator given by 
$$
(h_1,h_2)_c=\tilde h_1\tilde h_2\tilde h_1^{-1}\tilde h_2^{-1},
$$
where $\tilde h_i$ is any lifting of $h_i$ to $E$. If $H$ is abelian then the commutator $(h_1, h_2)_c$ depends only on the isomorphism class of the central extension. 

 Note that $T(\bF)=X_*(T)\otimes_{\ZZ} \bF^*$. For $f_i\in \bF^*$ and $\lambda_i\in X_*(T)$ the commutator for the central extension (\ref{ext_E_a}) is given by
$$
(\lambda_1\otimes f_1, \lambda_2\otimes f_2)_c=(f_1, f_2)_{st}^{2\check{h}(\lambda_1,\lambda_2)}
$$ 
Indeed, for $\lambda\in X_*(T)$, $f\in \bF^*$ the fibre of $\fL$ at $\lambda\otimes f$ identifies as $\ZZ/2\ZZ$-graded line with
$$
\otimes_{\check{\alpha}\in R^*}\det(\fg^{\check{\alpha}}(\bO): f^{\langle\lambda, \check{\alpha}\rangle}\fg^{\check{\alpha}}(\bO))
$$
and, using Lemma~\ref{lem_Coxeter}, it suffices to calculate the commutator of the canonical central extension of $\Gm(\bF)$ by $\Gm$ given by the relative determinant. But the latter commutator is given by the tame symbol (cf. \cite{DeB}, 12.13, p. 82).
We learn that the commutator for the central extension (\ref{ext_Mumford_sc_case}) is given on the torus by
$$
(\lambda_1\otimes f_1, \lambda_2\otimes f_2)_c=(f_1, f_2)_{st}^{(\lambda_1,\lambda_2)}
$$ 
for $\lambda_i\in X_*(T^{sc}), f_i\in\bF^*$. 

 We will check that the commutators corresponding to (\ref{ext_Mumford_sc_case}) and to (\ref{ext_sc_case_k_points}) are the same on $T^{sc}(\bF)$. The commutator for (\ref{ext_sc_case_k_points}) can by calculated using, for example, (\cite{DeB}, Proposition~11.11, p. 77). Namely, consider first the case of $G^{sc}=\SL_2$. In this case identify $T^{sc}$ with $\Gm$ via the positive coroot $\alpha: \Gm\,\iso\, T^{sc}$ then the commutator for (\ref{ext_sc_case_k_points}) becomes
$$
(f_1, f_2)_c=(f_1, f_2)^2_{st}
$$
Indeed, for $h_i\in T^{sc}$ consider in the notation of (\cite{DeB}, formula (11.1.4) on p. 73) Steinberg's cocycle
$c(h_1,h_2)\in K_2$. The image of $c(h_1,h_2)$ under the tame symbol $K_2(\bF)\to \sk^*$ equals $(h_1, h_2)_{st}$. So, the commutator $(f_1, f_2)_c$ is the image of 
$$
\frac{c(f_1, f_2)}{c(f_2, f_1)}\in K_2(\bF)
$$
under the tame symbol $K_2(\bF)\to \sk^*$. For $G^{sc}=\SL_2$ our assertion follows.

 The general case can be reduced to the case $G^{sc}=\SL_2$ by restricting to the $\SL_2$-subgroups $S_{\check{\alpha}}\subset G$ corresponding to the roots $\check{\alpha}$ as in (\cite{DeB}, Section~11.2, p.~74). 
We are done.  

\begin{rem}
\label{rem_parity_iota}
 For $\lambda\in X_*(T^{sc})$ we have $(\lambda,\lambda)\in 2\ZZ$. Indeed, $(\lambda,\lambda)=Q(2\lambda)-2Q(\lambda)=2Q(\lambda)$.
\end{rem}

\subsection{Proof of Proposition~\ref{prop_1}} The idea of the argument below was communicated to us by Drinfeld. 

 For $m\in\ZZ$ write $\bar G_{mQ}$ for the $m$-th multiple of the central extension $\bar G_Q$. There is a canonical action $\delta_0$ of $G(\bF)$ on the exact sequence 
\begin{equation}
\label{three}
1\to \Gm\to \bar G_{mQ}\to G^{sc}(\bF)\to 1
\end{equation}
such that $G(\bF)$ acts trivially on $\Gm$ and by conjugation on $G^{sc}(\bF)$. Indeed, we know that the extension (\ref{ext_Mumford_sc_case}) comes from the canonical extension (\ref{ext_can_K2}), so that the automorphisms of $G^{sc}$ act on it. Write $\delta_0: G(\bF)\times \bar G_{mQ}\to \bar G_{mQ}$ for the action map.
 
  If $\lambda\in X_*(T^{sc}), \mu\in X_*(T)$, $f,g\in \bF^*$ then $\mu\otimes g\in T(\bF)\subset G(\bF)$ acts on the fibre of $\bar G_{mQ}$ over $\lambda\otimes f\in T^{sc}(\bF)$ via $\delta_0$ as a multiplication by
$$
(g,f)_{st}^{m(\mu,\lambda)}
$$  
This is a kind of `explanation' of the fact that the form $(.,.)$ initially defined on $X_*(T^{sc})$ extends by linearity to a form $(.,.): X_*(T)\times X_*(T^{sc})\to\ZZ$ taking values  in $\ZZ$ and not just in $\BQ$.
 
\subsubsection{}  The isomorphism classes of central extensions 
\begin{equation}
\label{two}
1\to \Gm\to \cE_b\to T(\bF)\to 1
\end{equation}
are classified by symmetric bilinear forms $(.,.)_b: X_*(T)\times X_*(T)\to\ZZ$, namely for the corresponding extension (\ref{two}) we have
$$
(\lambda_1\otimes f_1, \lambda_2\otimes f_2)_c=(f_1,f_2)_{st}^{(\lambda_1,\lambda_2)_b}
$$ 
for $f_i\in\bF^*$, $\lambda_i\in X_*(T)$.  
 
 The group of automorphisms of the central extension (\ref{two}) is $\Hom(T(\bF), \Gm)$. Since $T(\bF)$ is abelian, the commutator $(\lambda_1\otimes f_1, \lambda_2\otimes f_2)_c$ is invariant under these automorphisms. The extension (\ref{two}) admits a (non unique) splitting over $T(\bO)$. 
 
\subsubsection{} The group $T$ acts on $G^{sc}$ by conjugation, let $\tilde G=G^{sc}\ltimes T$ denote the corresponding semi-direct product. The map $G^{sc}\ltimes T\to G$, $(g,t)\mapsto \bar g t$, where $\bar g$ is the image of $g\in G^{sc}$ in $G$, yields an exact sequence $1\to T^{sc}\to \tilde G\to G\to 1$. Hence, an exact sequence
\begin{equation}
\label{seq_from_Drinfeld}
1\to T^{sc}(\bF)\to \tilde G(\bF)\to G(\bF)\to 1
\end{equation}

The category of central extensions of $G(\bF)$ by $\Gm$ is equivalent to the category of pairs: a central extension 
\begin{equation}
\label{one}
1\to \Gm\to ?\to \tilde G(\bF)\to 1
\end{equation}
together with a splitting of its pull-back under $T^{sc}(\bF)\to \tilde G(\bF)$.  

 By a slight (we have to drop off the assumption of being of finite type) generalization of (\cite{DeB}, Construction~1.7), the category of central  
extensions (\ref{one}) is equivalent to the category of triples: central extensions (\ref{three})
and (\ref{two}) together with an action $\delta$ of $T(\bF)$ on $\bar G_{mQ}$ extending the action of $T(\bF)$ on $G^{sc}(F)$ by conjugation. Write $\delta: T(\bF)\times \bar G_{mQ}\to \bar G_{mQ}$ for the action map.

 Since (\ref{three}) has no automorphisms, $\delta$ coincides with the restriction of $\delta_0$ to $T(\bF)\times \bar G_{mQ}$. Write 
\begin{equation}
\label{ext_restriction_Mumford_to_T} 
1\to \Gm\to \bar G_{mQ}^T\to T^{sc}(\bF)\to 1
\end{equation}
for the restriction of (\ref{three}) to $T^{sc}(\bF)$. 
We conclude that the category of central extensions of $G(\bF)$ by $\Gm$ is equivalent to the category of pairs: a central extension (\ref{two}) together with an isomorphism of its restriction to $T^{sc}(\bF)$ with (\ref{ext_restriction_Mumford_to_T}). Clearly, the corresponding form $(.,.)_b: X_*(T)\times X_*(T)\to \ZZ$ is given by $(\lambda_1,\lambda_2)_b=m(\lambda_1,\lambda_2)$. Proposition~\ref{prop_1} follows.
 
\medskip

\section{Proof of Theorem~\ref{thm_main}}
\label{section_proof}

\subsection{Functors $F'_P$}
\label{subsection_functors_F_P}  Let $P$ be a parabolic subgroup of $G$ containing $B$. Let $M\subset P$ be a Levi subgroup containing $T$. Write
$$
1\to \Gm\to E^a_M\to M(\bF)\to 1
$$
for the restriction of (\ref{ext_E_a}) to $M(\bF)$, it is equipped with an action of $\Aut^0(\bO)$ and a section over $M(\bO)$ coming from the corresponding objects for (\ref{ext_E_a}). 
Let $\Perv_{M,G,N}$ denote the category of $M(\bO)$-equivariant $\Qlb$-perverse sheaves on $E^a_M/M(\bO)$ with $\Gm$-monodromy $\zeta_a$. Set
$$
\PPerv_{M,G,N}=\Perv_{M,G,N}[-1]\subset \D(E^a_M/M(\bO))
$$
 Write $\fL_{M,G}$ for the restriction of $\fL$ under $\Gr_M\to\Gr_G$, we equip it with the action of $\Aut^0(\bO)$ coming from that on $\fL$.  
 
 Write $\Gr_M$ for the affine grassmanian for $M$. The connected components of $\Gr_M$ are indexed by $\pi_1(M)$. For $\theta\in\pi_1(M)$ write $\Gr_M^{\theta}$ for the connected component of $\Gr_M$ containing $t^{\lambda}M(\bO)$ for any coweight $\lambda$ whose image in $\pi_1(M)$ is $\theta$. The diagram $M\gets P\hook{} G$ yields the following diagram of affine grassmanians 
$$
\Gr_M\getsup{\ft_P} \Gr_P\toup{\fs_P} \Gr_G
$$
The map $\ft_P$ yields a bijection between the connected components of $\Gr_P$ and those of $\Gr_M$. Let $\Gr_P^{\theta}$ be the connected component of $\Gr_P$ such that $\ft_P$ restricts to a map $\ft_P^{\theta}: \Gr_P^{\theta}\to \Gr_M^{\theta}$. Write $\fs^{\theta}_P: \Gr_P^{\theta}\to\Gr_G$ for the restriction of $\fs_P$. The restriction of $\fs^{\theta}_P$ to $(\Gr^{\theta}_P)_{red}$ is a closed immersion. 

 The section $M\hook{}P$ yields a section $\fr_P: \Gr_M\to\Gr_P$ of $\ft_P$. By abuse of notation, we write 
$$
\Gra_M\toup{\fr_P}\Gra_P\toup{\fs_P}\Gra_G
$$ 
for the diagram obtained from $\Gr_M\toup{\fr_P}\Gr_P
 \toup{\fs_P}\Gr_G$ by the base change $\Gra_G\to\Gr_G$. Clearly, $\ft_P$ lifts naturally to a map $\ft_P: \Gra_P\to\Gra_M$. Define the functor 
$$
F'_P: \PPerv_{G,N}\to \D(\Gra_M)
$$ 
by
$F'_P(K)=\ft_{P !}\fs_P^*K$. 
Write $\Gra_M^{\theta}$ (resp., $\Gra_P^{\theta}$) for the connected component of $\Gra_M$ (resp., $\Gra_P$) over 
$\Gr_M^{\theta}$ (resp., $\Gr_P^{\theta}$). Write 
$$
\PPerv^{\theta}_{M,G,N}\subset\PPerv_{M,G,N}
$$ 
for the full subcategory of objects that vanish off $\Gra_M^{\theta}$. Set 
$$
\PPerv'_{M,G,N}=\mathop{\oplus}\limits_{\theta\in\pi_1(M)} \PPerv^{\theta}_{M,G,N}[\langle\theta, 2\check{\rho}_M-2\check{\rho}\rangle]
$$ 
As in (\cite{BD}, 5.3.29) one shows that $F'_P$ sends $\PPerv_{G,N}$ to $\PPerv'_{M,G,N}$ (cf. also \cite{L}, appendix A.4). 

 The above construction applied to the Borel subgroup yields a functor $F'_B: \PPerv_{G,N}\to \PPerv'_{T,G,N}$. 

 Let $B(M)\subset M$ be a Borel subgroup containing $T$ such that the preimage of $B(M)$ under $P\to M$ equals $B$. The functor $F'_{B(M)}: \PPerv'_{M,G,N}\to \D(\Gra_T)$ is defined as follows. As above, the inclusions $T\subset B(M)\subset M$ yield a diagram 
\begin{equation}
\label{diag_Gr_for_Levi}
\Gr_T\toup{\fr_{B(M)}}\Gr_{B(M)} \toup{\fs_{B(M)}} \Gr_M
\end{equation}
Write 
$$
\Gra_T\toup{\fr_{B(M)}} \Gra_{B(M)}\toup{\fs_{B(M)}} \Gra_M
$$
for the diagram obtained from (\ref{diag_Gr_for_Levi}) by the base change $\Gra_M\to\Gr_M$. The projection $B(M)\to T$ yields $\ft_{B(M)}: \Gr_{B(M)}\to\Gr_T$, which lifts naturally to $\ft_{B(M)}: \Gra_{B(M)}\to\Gra_T$. For $K\in \PPerv'_{M,G,N}$ set
$$
F'_{B(M)}(K)=(\ft_{B(M)})_!\fs_{B(M)}^*K
$$
As in (\cite{BD}, 5.3.29), one shows that $F'_{B(M)}$ is a functor 
$$
F'_{B(M)}: \PPerv'_{M,G,N}\to \PPerv'_{T,G,N}
$$ 
By base change, we have canonically
\begin{equation}
\label{iso_tens_functors_for_P}
F'_{B(M)}\comp F'_P\,\iso\, F'_B
\end{equation} 

\subsubsection{} Write $X_*^{+M}(T)\subset X_*(T)$ for the coweights of $T$ dominant for $M$. For $\lambda\in X_*^{+M}(T)$ denote by $\Gr_M^{\lambda}$ the $M(\bO)$-orbit through $t^{\lambda}M(\bO)$. Let $\Gra_M^{\lambda}$ be the preimage of $\Gr_M^{\lambda}$ under $\Gra_M\to\Gr_M$. The $M$-orbit on $\Gr_M$ through $t^{\lambda}M(\bO)$ is isomorphic to a partial flag variety $\CB_M^{\lambda}=M/P^{\lambda}_M$, where the Levi subgroup of $P^{\lambda}_M$ has the Weyl group 
coinciding with the stabilizer of $\lambda$ in $W_M$. Write $\tilde\omega_{M,\lambda}: \Gr^{\lambda}_M\to\CB^{\lambda}_M$ for the projection. As in Lemma~\ref{lem_fibres_fL}, one gets a $\ZZ/2\ZZ$-graded isomorphism
$$
\fL_{M,G}\mid_{\Gr_M^{\lambda}}\,\iso\,\Omega_{\bar c}^{\check{h}(\lambda,\lambda)}\otimes \tilde\omega_{M,\lambda}^*\CO(2\check{h}\iota(\lambda)),
$$
here the line bundles $\CO(\check{\nu})$ on $\CB_M^{\lambda}$ are defined as in Section~\ref{sect_notations}.
Set
$$
X^{*+}_M(\check{T}_N)=\{\lambda\in  X_*^{+M}(T)\mid d\iota(\lambda)\in N X^*(T)\}
$$
As for $G$ itself, for $\lambda\in X_*^{+M}(T)$ the scheme $\Gra_M^{\lambda}$ admits a $M(\bO)$-equivariant local system with $\Gm$-monodromy $\zeta_a$ iff $\lambda\in X^{*+}_M(\check{T}_N)$. For $\lambda\in X^{*+}_M(\check{T}_N)$ denote by $\CA^{\lambda}_{M,\cE}$ the irreducible object of $\PPerv_{M,G,N}$ defined as in Section~\ref{sect_notations}.

\subsubsection{More tensor structures} One equips $\PPerv_{M,G,N}$ and $\PPerv'_{M,G,N}$ with a convolution product as in Section~\ref{sect_convolution}. Let us define the commutativity constraint on these categories via fusion. 

 Recall the line bundles $\fL_{X^m}$ on $\Gr_{G,X^m}$ from Section~\ref{sect_Fusion}. For the convenience of the reader we remind the \select{factorization structure} on these line bundles, which allowed to do fusion for $\PPerv_{G,N}$. 
 
  For a surjective map of finite sets $\alpha: J\to I$ one has a cartesian square
$$
\begin{array}{ccc}
\Gr_{G, X^{I}} & \toup{{\tilde\triangle}^{\alpha}} & 
\Gr_{G, X^J}\\
\downarrow && \downarrow\\
X^I & \toup{\triangle^{\alpha}} & X^J,
\end{array}
$$ 
where $\triangle^{\alpha}$ is the corresponding diagonal. We have canonically $(\tilde\triangle^{\alpha})^*\fL_{X^J}\,\iso\, \fL_{X^I}$. 

 Write $\nu^{\alpha}: U^{\alpha}\hook{} X^J$ for the open subscheme given by the condition that the divisors $D_i$ do not intersect pairwise, where $D_i=\sum_{j\in J, \;\alpha(j)=i} x_j$ for $(x_j)\in X^J$.  We have a cartesian square
$$
\begin{array}{ccc}
(\prod_{i\in I}\Gr_{G, X^{\alpha^{-1}(i)}})\mid_{U^{\alpha}} & \hook{{\tilde\nu}^{\alpha}} & \Gr_{G, X^J}\\
\downarrow && \downarrow\\
U^{\alpha} & \toup{\nu^{\alpha}} & X^J
\end{array}
$$
We have canonically 
$$
(\tilde\nu^{\alpha})^*\fL_{X^J}\,\iso\, (\boxtimes_{i\in I} \fL_{X^{\alpha^{-1}(i)}})\mid_{U^{\alpha}}
$$  
Let $\fL_{M,G,X^m}$ be the restriction of $\fL_{X^m}$ under the map $\Gr_{M, X^m}\to \Gr_{G,X^m}$ induced by $M\hook{}G$. The collection $\{\cL_{M,G, X^m}\}$ is endowed with the induced factorization structure. Write $\Gra_{M,G,X^m}$ for  the punctured line bundle of $\fL_{M,G,X^m}$. 

 Let $M_{X^m}$ be the group scheme over $X^m$ classifying  $\{(x_1,\ldots, x_m)\in X^m, \mu\}$, where $\mu$ is an automorphism of $\cF^0_M$ restricted to the formal neighborhood of $x_1\cup\ldots\cup x_m$ in $X$. The group scheme $M_{X^m}$ acts naturally on $\Gra_{M, G, X^m}$. Write $\Perv_{M,G,N,X^m}$ for the category of $M_{X^m}$-equivariant perverse sheaves on $\Gra_{M,G,X^m}$ with $\Gm$-monodromy $\zeta_a$. Set
$$
\PPerv_{M,G,N,X^m}=\Perv_{M,G,N,X^m}[-m-1]
$$ 

  Let $\Aut^0_2(\bO)$ act on $\Gra_M$ via its quotient $\Aut^0(\bO)$. Then every object of $\Perv_{M,G,N}$ admits a unique $\Aut_2^0(\bO)$-equivariant structure. Note that $\Gra_{M, X}\,\iso\, \hat X_2\times_{\Aut_2^0(\bO)}\Gra_M$. As in Section~\ref{sect_Fusion}, we get a fully faithful functor
$$
\tau^0: \PPerv_{M,G,N}\to\PPerv_{M,G,N,X}  
$$
Now we define the commutativity constraint on $\PPerv_{M,G,N}$ and $\PPerv'_{M,G,N}$ using the above factorization structure as in Section~\ref{sect_Fusion}. As in (\cite{BD}, 5.3.16) one checks that $\PPerv_{M,G,N}$ and $\PPerv'_{M,G,N}$ are symmetric monoidal categories. 

\begin{lem} The functors $F'_P$, $F'_{B(M)}$, $F'_B$ are tensor functors, and (\ref{iso_tens_functors_for_P}) is an isomorphism of tensor functors.
\end{lem}
\begin{Prf} We will only check that $F'_P$ is a tensor functor, the rest is similar.\\
1) Let us show that $F'_P$ is compatible with the convolution. Let $\Gr_{P, X^m}$ be the ind-scheme classifying $(x_1,\ldots,x_m)\in X^m$, a $P$-torsor $\cF_P$ on $X$, and a trivialization $\cF_P\,\iso\, \cF^0_P\mid_{X-\cup x_i}$. Write $\Gra_{P,X^m}$ for the ind-scheme obtained  from $\Gra_{G,X^m}$ by the base change $\Gr_{P,X^m}\to\Gr_{G,X^m}$. As in Section~\ref{subsection_functors_F_P}, we get a diagram 
$$
\Gra_{M,X^m} \getsup{\ft_{P,X^m}} \Gra_{P, X^m} \toup{\fs_{P, X^m}} \Gra_{G, X^m}
$$
and a functor 
$$
F'_{P,X^m}: \D(\Gra_{G, X^m})\to \D(\Gra_{M, G, X^m})
$$ 
given by $F'_{P,X^m}(K)=(\ft_{P,X^m})_!\fs_{P,X^m}^*$. For $i=1,2$ let $F_i\in\PPerv_{G,N}$ and $K_i=\tau^0 F_i$. 
Recall that $U\subset X^2$ is the complement to the diagonal. We have a natural diagram, where both squares are cartesian
$$
\begin{array}{ccc}
(\Gra_{G,X}\times\Gra_{G,X})\mid_U & \toup{\nu_{G,U}} & \Gra_{G,X^2}\\
\uparrow && \uparrow\lefteqn{\scriptstyle \fs_{P,X^2}}\\
(\Gra_{P,X}\times\Gra_{P,X})\mid_U & \to & \Gra_{P,X^2}\\
\downarrow && \downarrow\lefteqn{\scriptstyle \ft_{P, X^2}}\\
(\Gra_{M,X}\times\Gra_{M,X})\mid_U & \toup{\nu_{M,U}}  & \Gra_{M,X^2},
\end{array}
$$
and the maps $\nu_{G,U}$ and $\nu_{M,U}$ come from the above factorization structures. As in (\cite{L}, Proposition~14), one shows that $F'_{P,X^2}(K_1\ast_X K_2)$ is the Goresky-MacPherson extension from $\Gra_{M,X^2}\mid_U$. Now the isomorphism 
$$
\nu_{M,U}^*F'_{P,X^2}(K_1\ast_X K_2)\,\iso\, \tau^0(F'_P(F_1))\boxtimes \tau^0(F'_P(F_2))
$$ 
yields an isomorphism
$$
\epsilon_{12}: F'_{P,X^2}(K_1\ast_X K_2)\,\iso\, \tau^0(F'_P(F_1))\ast_X \tau^0(F'_P(F_2))
$$
Restricting it to the diagonal in $X$ one gets 
$$
\tau^0 (F'_P(F_1\ast F_2))\,\iso\,
\tau^0(F'_P(F_1)\ast F'_P(F_2))
$$
2) Let us check that $F'_P$ is compatible with the commutativity constraints. Recall that $\sigma$ is the involution of $X^2$ permuting the two coordinates. One has a commutative diagram, where the vertical arrows are canonical isomorphisms
$$
\begin{array}{ccc}
\sigma^* F'_{P,X^2}(K_1\ast_X K_2) & \toup{\sigma^*\epsilon_{12}} & \sigma^*(\tau^0(F'_P(F_1))\ast_X \tau^0(F'_P(F_2)))\\
\downarrow && \downarrow\\
F'_{P,X^2}(K_2\ast_X K_1) & \toup{\epsilon_{21}} & \tau^0(F'_P(F_2))\ast_X \tau^0(F'_P(F_1))
\end{array}
$$
Restricting it to the diagonal, one gets the desired compatibility isomorphism. 
\end{Prf}
  
\subsection{Fibre functor} 
\label{subsection_fibre_functor} In Section~\ref{subsection_functors_F_P} we introduced the ($\ZZ/2\ZZ$-graded purely of parity zero)
line bundle $\fL_{T,G}$ on $\Gr_T$. The action of $T(\bO)$ on this line bundle comes from the action of $T(\bO)$ on $E^a_T$ by left multpilication. 
The fibre of $\fL_{T,G}$ at $\nu\in X_*(T)$ is $T(\bO)$-equivariantly isomorphic to $\Omega_{\bar c}^{\check{h}\langle\nu, \iota(\nu)\rangle}$, where $T(\bO)$ acts on the latter space via $T(\bO)\to T\toup{2\check{h} \iota(\nu)}\Gm$.
Recall the torus $\check{T}_N$ introduced in Section~\ref{sect_main_result}. For $\nu\in X_*(T)$ the orbit $\Gra^{\nu}_T$ supports a nonzero object of $\PPerv_{T,G,N}$ iff $\nu\in X^*(\check{T}_N)$.

 For $\nu\in X^*(\check{T}_N)$ consider the map $a_{\nu}: \cE_{\bar c}^{d(\nu,\nu)/N}-\{0\}\to \Omega_{\bar c}^{\check{h}(\nu,\nu)}-\{0\}$ sending $x$ to $x^{2\check{h}N/d}$. For $K\in \PPerv_{T,G,N}$ the complex $a_{\nu}^*K$ is a constant sheaf placed in degree zero, so we view it as a vector space denoted $F_T^{\nu}(K)$. Then 
$$
F_T=\mathop{\oplus}\limits_{\nu\in X^*(\check{T}_N)} F_T^{\nu}
$$ 
is a fibre functor $\PPerv_{T,G,N}\to\Vect$. Let $\check{T}_N=\Spec k[X^*(\check{T}_N)]$ be the torus whose weight lattice is $X^*(\check{T}_N)$. By (\cite{DM}, Theorem 2.11), we get 
$$
\PPerv_{T,G,N}\,\iso\, \Rep(\check{T}_N)
$$
 
 For $\nu\in X^*(\check{T}_N)$ write $F_{B(M)}^{'\nu}$ for the functor $F'_{B(M)}$ followed by restriction to $\Gra_T^{\nu}$. Write $F^{\nu}_M: \PPerv_{M,G,N}\to\Vect$ for the functor 
$$
F^{\nu}_T  F_{B(M)}^{'\nu} [\langle\nu, 2\check{\rho}_M\rangle]
$$ 
For $\nu\in X^*(\check{T}_N)$ any $x\in \cE_{\bar c}^{d(\nu,\nu)/N}$ yields a section $a_{B(M),\nu}: \Gr_{B(M)}^{\nu}\to \Gra_{B(M)}^{\nu}$ of the projection $\Gra_{B(M)}^{\nu}\to\Gr^{\nu}_{B(M)}$ sending $x$ to $x^{2\check{h}N/d}$.

\begin{lem}  
\label{lem_F_nu(A_lambda)_has_base}
If $\nu\in X^*(\check{T}_N)$, $\lambda\in X^{*+}_M(\check{T}_N)$ then $F^{\nu}_M(\CA^{\lambda}_{M,\cE})$ has a canonical base consisting of those connected components of $$
\Gr_{B(M)}^{\nu}\cap \Gr^{\lambda}_M
$$ 
over which the (shifted) local system $a_{B(M),\nu}^*\CA^{\lambda}_{M,\cE}$ is constant. In particular, 
$F^{\lambda}_M(\CA^{\lambda}_{\cE})\,\iso\,\Qlb$ 
canonically, and $F^{w(\lambda)}(\CA^{\lambda}_{\cE})\,\iso\,\Qlb$ for $w\in W_M$.
\end{lem}
\begin{Prf}
The first claim is proved as in (\cite{MV}, Proposition~3.10). If $\lambda\in X^{*+}_M(\check{T}_N)$ then $\Gr^{w(\lambda)}_{B(M)}\cap \Gr^{\lambda}_M$ is an affine space (\cite{MV}, proof of Theorem~3.2), and the local system $a_{B(M),\lambda}^*\CA^{\lambda}_{M,\cE}$ is constant. 
\end{Prf}

\medskip

 Consider the following $\ZZ/2\ZZ$-grading on $\PPerv'_{M,G,N}$. For $\theta\in \pi_1(M)$ call an object of 
$
\PPerv^{\theta}_{M,G,N}[\langle\theta, 2\check{\rho}_M-2\check{\rho}\rangle]
$ 
\select{of parity} $\langle\theta, 2\check{\rho}\rangle\!\!\mod 2$. Write $E^{a,\theta}_M$ for the connected component of $E^a_M$ such that 
$$
E^{a,\theta}_M/M(\bO)=\Gra_M^{\theta}
$$ 
The product in $E^a_M$ is compatible with this $\ZZ/2\ZZ$-grading of $\pi_1(M)$, so the $\ZZ/2\ZZ$-grading we get on $\PPerv'_{M,G,N}$ is compatible with the tensor structure.  
In particular, for $M=G$ we get a $\ZZ/2\ZZ$-grading on $\PPerv_{G,N}$. If $(\Gra_P^{\theta})_{red}$ is contained in the connected component $\Gra_G^{\bar\theta}$ of $\Gra_G$ then $\bar\theta$ is the image of $\theta$ in $\pi_1(G)$. So, the functors $F'_P$ and $F'_{B(M)}$ are compatible with these gradings.

 Write $\Vect^{\epsilon}$ for the tensor category of $\ZZ/2\ZZ$-graded vector spaces. Let $\PPerv_{M,G,N}^{\natural}$ be the category of even objects in $\PPerv'_{M,G,N}\otimes\Vect^{\epsilon}$. Let $\PPerv_{G,N}^{\natural}$ be the category of even objects in $\PPerv_{G,N}\otimes\Vect^{\epsilon}$. We get a canonical equivalence of tensor categories
$
sh: \PPerv^{\natural}_{T,G,N}\,\iso\,\PPerv_{T,G,N}
$
The functors $F'_{B(M)}, F'_P, F'_B$ yield tensor functors 
\begin{equation}
\label{diag_tens_natural_three}
\PPerv^{\natural}_{G,N}\,\toup{F^{\natural}_P}\, \PPerv_{M,G,N}^{\natural}\,\toup{F^{\natural}_{B(M)}} \,\PPerv_{T,G,N}^{\natural}
\end{equation}
whose composition is $F^{\natural}_B$.
Write $F^{\natural}: \PPerv_{G,N}^{\natural}\to\Vect$ for $F_T\comp sh\comp F^{\natural}_B$. By Lemma~\ref{lem_F_nu(A_lambda)_has_base}, $F^{\natural}$ does not annihilate a nonzero object, so it is faithful. By Remark~\ref{rem_rigidity}, $\PPerv^{\natural}_{G,N}$ is a rigid abelian tensor category, so $F^{\natural}$ is a fibre functor. By (\cite{DM}, Theorem~2.11), $\Aut^{\otimes}(F^{\natural})$ is represented by an affine group scheme $\check{G}_N$, and we have an equivalence of tensor categories 
\begin{equation}
\label{equiv_PPerv_natural_Rep_smth}  
  \PPerv^{\natural}_{G,N}\,\iso\,\Rep(\check{G}_N)
\end{equation}

 An analog of Remark~\ref{rem_rigidity} holds also for $M$, so $F_T\comp sh\comp F^{\natural}_{B(M)}: \PPerv^{\natural}_{M,G,N}\to \Vect$ is a fibre functor that yields an affine group scheme $\check{M}_N$ and an equivalence of tensor categories $\PPerv^{\natural}_{M,G,N}\,\iso\, \Rep(\check{M}_N)$. The diagram (\ref{diag_tens_natural_three}) yields homomorphisms $\check{T}_N\to\check{M}_N\to\check{G}_N$. Since $X^{*+}(\check{T}_N)$ does not contain a nontrivial subgroup, $\check{G}_N$ is semisimple of rank equal to the rank of $G$.
  
\subsection{Structure of $\check{G}_N$}   
   
\begin{lem}
\label{lem_lambda_plus_mu}
 If $\lambda,\mu\in X^{*+}(\check{T}_N)$ then $\CA^{\lambda+\mu}_{\cE}$ appears in $\CA^{\lambda}_{\cE}\ast \CA^{\mu}_{\cE}$ with multiplicity one.
\end{lem}
\begin{Prf} Write $\bar E^{a,\lambda}$ (resp., $E^{a,\lambda}$) for the preimage of $\ov{\Gra}^{\lambda}_G$ (resp., $\Gra^{\lambda}_G$) under $E^a\to\Gra_G$, $x\mapsto xG(\bO)$. Write $m^{\lambda,\mu}:\bar E^{a,\lambda}\times_{G(\bO)\times\Gm} \ov{\Gra}^\mu_G\to \ov{\Gra}^{\lambda+\mu}_G$ for the convolution diagram as in Section~\ref{sect_convolution}. If $W$ is the preimage of $\Gra_G^{\lambda+\mu}$ under $m^{\lambda,\mu}$ then $m^{\lambda,\mu}: W\to \Gra_G^{\lambda+\mu}$ is an isomorphism, and $W$ is open in $E^{a,\lambda}\times_{G(\bO)\times\Gm} \Gra^{\mu}_G$.  
\end{Prf}
   
\medskip

   Clearly, if $X^{*+}(\check{T}_N)$ is a $\ZZ_+$-span of
$\lambda_1,\ldots,\lambda_r$ then $\oplus_i\CA^{\lambda_i}_{\CE}$ is a tensor generator for $\PPerv_{G,N}$. So, $\check{G}_N$ is algebraic by (\cite{DM}, Proposition~2.20). 
By Lemma~\ref{lem_F_nu(A_lambda)_has_base}, for $\mu\in X^{*+}(\check{T}_N)$ and $w\in W$ the weight $w(\mu)$ of $\check{T}_N$ appears in $F^{\natural}(\CA^{\mu}_{\CE})$. So, $\check{T}_N$ is closed in $\check{G}_N$ by (\cite{DM}, Proposition~2.21). By Lemma~\ref{lem_lambda_plus_mu}, there is no tensor subcategory of $\PPerv_{G,N}$ whose objects are direct sums of finitely many fixed irreducible objects, so $\check{G}_N$ is connected by (\cite{DM}, Corollary~2.22). Since $\PPerv_{G,N}$ is semisimple, $\check{G}_N$ is reductive by (\cite{DM}, Proposition~2.23). We will use the following.
 
\begin{lem} 
\label{lem_general_Borel}
Let $\GG$ be a connected reductive group with a maximal torus $\TT\subset\GG$. Let $\check{\Lambda}^+$ be a subsemigroup in the group $\check{\Lambda}$ of weights of $\TT$. Assume that we are given a bijection $\nu\mapsto V^{\nu}$ between $\check{\Lambda}^+$ and the set of irreducible representations of $\GG$ such that the following holds:
\begin{itemize}
\item if $\nu\in\check{\Lambda}^+$ then the $\nu$-weight space $L^{\nu}$ of $\TT$ in $V^{\nu}$ is of dimension one;
\item if $\nu_1,\nu_2\in\check{\Lambda}^+$ then $V^{\nu_1+\nu_2}$ appears with multiplicity one in $V^{\nu_1}\otimes V^{\nu_2}$, and the subspace $L^{\nu_1}\otimes L^{\nu_2}\subset V^{\nu_1}\otimes V^{\nu_2}$ coincides with the image of $L^{\nu_1+\nu_2}\hook{} V^{\nu_1+\nu_2}\hook{} V^{\nu_1}\otimes V^{\nu_2}$.
\end{itemize}
Then there is a unique Borel subgroup $\TT\subset\BB\subset \GG$ such that $\check{\Lambda}^+$ is the set of dominant weights for $\BB$. $\square$
\end{lem}
 
Write $V^{\nu}$ for the irreducible representation of $\check{G}_N$ corresponding to $\CA^{\nu}_{\cE}$ via (\ref{equiv_PPerv_natural_Rep_smth}). 

\begin{lem}
\label{lem_Borel_exists}
 The torus $\check{T}_N$ is maximal in $\check{G}_N$. There is a unique Borel subgroup $\check{T}_N \subset \check{B}_N\subset \check{G}_N$ whose set of dominant weights coincides with $X^{*+}(\check{T}_N)$.
\end{lem}
\begin{Prf}
Let $T'\subset \check{G}_N$ be a maximal torus containing $\check{T}_N$. By Lemma~\ref{lem_F_nu(A_lambda)_has_base}, for each $\nu\in X^{*+}(\check{T}_N)$ there is a unique character $\nu'$ of $T'$ such that the composition $\check{T}_N\to T'\toup{\nu'}\Gm$ is $\nu$, and the $T'$-weight $\nu'$ appears in $V^{\nu}$. Clearly, $\nu\mapsto \nu'$ is a homomorphism of semigroups, and we can apply Lemma~\ref{lem_general_Borel}. Since $\nu\mapsto\nu'$ is a bijection between $X^{*+}(\check{T}_N)$ and the dominant weights of $\check{B}_N$, $\check{T}_N$ is maximal.
\end{Prf}   

\medskip

 Applying similar arguments for $M$, one checks that $\check{M}_N$ is reductive, and $\check{T}_N\to\check{M}_N\to\check{G}_N$ are closed immersions, so $\check{M}_N$ is a Levi subgroup of $\check{G}_N$. 
 
\subsection{Rank one} Let $M$ be the subminimal Levi subgroup of $G$ corresponding to the simple root $\check{\alpha}_i$. As in Lemma~\ref{lem_Borel_exists}, there is a unique Borel subgroup $\check{T}_N\subset \check{B}(M)_N\subset \check{M}_N$ whose set of dominant weights is $X^{*+}_M(\check{T}_N)$. View $\check{\alpha}_i$ as a coweight of $\check{T}_N$. Then 
$$
\{\check{\nu}\in  X_*(\check{T}_N)\mid \langle\lambda,\check{\nu}\rangle\ge 0\;\mbox{for all}\; \lambda\in X^{*+}_M(\check{T}_N)\}
$$ 
is a $\ZZ_+$-span of a multiple of $\check{\alpha}_i$. So, $\check{M}_N$ is of semisimple rank one, and its unique simple coroot is of the form $\check{\alpha}_i/\kappa_i$ for some $\kappa_i\in\QQ$, $\kappa_i>0$. 

 Take any $\lambda\in X^{*+}_M(\check{T}_N)$ with $\langle\lambda, \check{\alpha}_i\rangle >0$. Write $s_i\in W$ for the simple reflection corresponding to $\check{\alpha}_i$. By Lemma~\ref{lem_F_nu(A_lambda)_has_base}, $F^{\lambda}_M(\CA^{\lambda}_{M,\CE})$ and $F^{s_i(\lambda)}_M(\CA^{\lambda}_{M,\CE})$ do not vanish, so $\lambda-s_i(\lambda)$ is a multiple of the positive root of $\check{M}_N$. So, this positive root is $\kappa_i\alpha_i$. We learn that the simple reflection for $(\check{M}_N, \check{T}_N)$ acts on $X^*(\check{T}_N)$ as $s_i$. So, the Weyl groups of $G$ and  of $\check{G}_N$, viewed as subgroups of $\Aut(X^*(\check{T}_N))$ are the same. We must show that $\kappa_i=\delta_i$.
 
  Recall that the scheme $\Gr_{B(M)}^{\nu}\cap \Gr^{\lambda}_M$ is non empty iff 
$$
\nu=\lambda, \lambda-\alpha_i,\lambda-2\alpha_i,\ldots, \lambda-\langle\lambda,\check{\alpha}_i\rangle\alpha_i
$$ 
For $0<k<  \langle\lambda,\check{\alpha}_i\rangle$ and $\nu=\lambda-k\alpha_i$ one has 
$$
\Gr_{B(M)}^{\nu}\cap \Gr^{\lambda}_M\,\iso\,\Gm\times \AA^{\langle\lambda,\check{\alpha}_i\rangle-k-1}
$$
    
   Write $M_0$ for the derived group of $M$, let $T_0\subset M_0$ be the maximal torus such that $T_0\subset T$. Consider the central extension 
$$
1\to\Gm\to E^a_{M_0}\to M_0(\bF)\to 1
$$ 
obtained by pulling back of $1\to \Gm\to E^a_M\to M(\bF)\to 1$ via $M_0(\bF)\to M(\bF)$. It corresponds to the restriction of the bilinear form $2\check{h}\iota$ under 
$$
X_*(T_0)\times X_*(T_0)\subset X_*(T)\times X_*(T)
$$ 
So, $E^a_{M_0}/M_0(\bO)\to \Gr_{M_0}$ is isomorphic to the punctured total space of 
$\CL_{M_0}^{\check{h}(\alpha_i,\alpha_i)}$, where $\CL_{M_0}$ is an ample generator of the Picard group of (each connected component of) $\Gr_{M_0}$. 
 
 Assume that $\lambda=a\alpha_i$ with $a>0, a\in\ZZ$ such that $\lambda\in X^*(\check{T}_N)$. Let $\nu=b\alpha_i$ with $b\in\ZZ$ satisfy $-\lambda<\nu< \lambda$. 
 
 Write $U\subset M(\bF)$ for the 1-parameter unipotent subgroup corresponding to the affine root space $t^{-a+b}\fg_{\check{\alpha}_i}$. Let $Y$ be the closure of the $U$-orbit through $t^{\nu}M(\bO)$ in $\Gr_M$. It is a $T$-stable subscheme $Y\,\iso\,\PP^1$, the $T$-fixed points in $Y$ are
$t^{\nu}M(\bO)$ and $t^{-\lambda}M(\bO)$. The restriction of $\CL_{M_0}$ to $Y$ identifies with $\CO_{\PP^1}(a+b)$.
The section 
$$
a_{B(M),\nu}: \Gr^{\nu}_{B(M)}\to \Gra^{\nu}_{B(M)}
$$ 
viewed as a section of the line bundle $\CL_{M_0}^{\check{h}(\alpha_i,\alpha_i)}$ over $Y$ will vanish only at $t^{-\lambda}M(\bO)$ with multiplicity 
$(a+b)\check{h}(\alpha_i,\alpha_i)$. 
So, the local system $a_{B(M),\nu}^*\CA^{\lambda}_{M,\CE}$ will have the $\Gm$-monodromy $\zeta_a^{(a+b)\check{h}(\alpha_i,\alpha_i)}$. This local system is trivial iff 
$$
(a+b)\check{h}(\alpha_i,\alpha_i)\in \frac{2\check{h}N}{d}\ZZ
$$ We may assume that $\frac{da}{2N}(\alpha_i,\alpha_i)\in\ZZ$. Then the above condition is equivalent to $b\in \frac{2N}{d(\alpha_i,\alpha_i)}\ZZ$.  
The smallest positive integer $b$ satisfying this condition is $\delta_i$. So, $\kappa_i=\delta_i$. Theorem~\ref{thm_main} is proved.

\bigskip

\footnotesize{
{\bf M.F.}: IMU, IITP, and State University Higher School of Economy,\\
Department of Mathematics, 20 Myasnitskaya st, Moscow 101000 Russia;\\
%\hphantom{x}\ab\, 
{\tt fnklberg@gmail.com}}

\footnotesize{
{\bf S.L.}: Institut \'Elie Cartan Nancy (Math\'ematiques), Universit\'e Henri Poincar\'e Nancy 1, 
  B.P. 239, F-54506 Vandoeuvre-l\`es-Nancy Cedex, France\\
{\tt Sergey.Lysenko@iecn.u-nancy.fr}}

\end{document}